\documentclass[a4paper,fleqn,10pt,twoside,draft]{article}

\usepackage[utf8]{inputenc}

\input{layout_commands}
\usepackage{amsmath}
\usepackage{amssymb}
\usepackage{accents}

\newcounter{mynumcounter}
	       {%
		 \begin{list}{(\roman{mynumcounter})\hspace*{\fill}}%
		   {
		     \setlength{\topsep}{0cm}
		     \setlength{\partopsep}{0cm}
		     \setlength{\itemsep}{0ex}
		     \setlength{\parsep}{0cm}
		     \setlength{\leftmargin}{0cm}
		     \setlength{\itemindent}{8mm}		   
		     \setlength{\labelsep}{3mm}
		     \setlength{\labelwidth}{5mm}
		     \usecounter{mynumcounter}
		   }%
	       }
	       {\end{list}}

\newcommand{\eps}{\varepsilon}
\newcommand{\del}{\partial}

\newcommand{\dd}[2]{\frac{\del #1}{\del #2}}
\newcommand{\ddeval}[3]{\left.\dd{#1}{#2}\right|_{#3}}

\newcommand{\DD}[2]{\frac{d #1}{d #2}}
\newcommand{\DDeval}[3]{\left.\DD{#1}{#2}\right|_{#3}}

\newcommand{\Ric}{\operatorname{Ric}}
\newcommand{\Scal}{\operatorname{Scal}}

\newcommand{\diam}{\operatorname{diam}}
\renewcommand{\div}{\operatorname{div}}

\newcommand{\IC}{\mathbf{C}}

\newcommand{\IR}{\mathbf{R}}

\newcommand{\CA}{\mathcal{A}}
\newcommand{\CB}{\mathcal{B}}

\newcommand{\CV}{\mathcal{V}}
\newcommand{\CW}{\mathcal{W}}
\newcommand{\CX}{\mathcal{X}}

\newcommand{\al}{\alpha}
\newcommand{\be}{\beta}
\newcommand{\ga}{\gamma}
\newcommand{\de}{\delta}

\newcommand{\Id}{\operatorname{Id}}

\newcommand{\rmd}{\mathrm{d}}

\newcommand{\dmu}{\,\rmd\mu}
\newcommand{\dx}{\,\rmd x}

\newcommand{\ra}{\rangle}
\newcommand{\la}{\langle}
\newcommand{\inj}{\mathrm{inj}}

\newcommand{\Connection}[1]{\smash{\sideset{^{#1}}{}{\mathop\nabla\nolimits}}}

\newcommand{\nabM}{\!\Connection{M}}

\newcommand{\Riemann}[1]{\smash{\sideset{^{#1}}{}{\mathop\mathrm{Rm}\nolimits}}}

\newcommand{\Riem}{\!\Riemann{}}

\newcommand{\Ricci}[1]{\smash{\sideset{^{#1}}{}{\mathop\mathrm{Rc}\nolimits}}}

\newcommand{\RicSig}{\!\Ricci{\Sigma}}

\newcommand{\Scalarcurv}[1]{\smash{\sideset{^{#1}}{}{\mathop\mathrm{Sc}\nolimits}}}

\newcommand{\ScalSig}{\!\Scalarcurv{\Sigma}}

\newcommand{\Acirc}{\accentset{\circ}{A}}

\newcommand{\half}{\tfrac{1}{2}}

\input{draft_commands}

\newcommand{\dsig}{\mathrm{d}\sigma}
\renewcommand{\Scal}{\operatorname{Sc}}
\newcommand{\weakto}{\rightharpoonup}
\newcommand{\weaktostar}{\stackrel{*}{\weakto}}
\renewcommand{\IR}{\mathbb{R}}

\renewcommand{\IC}{\mathbb{C}}

\title{\titlefamily\Huge Minimizers of the Willmore functional with a small area constraint}
\author{
\authname{Tobias Lamm}
  \authaddress
  {
   Goethe-Universit\"at Frankfurt, 
   Institut f\"ur Mathematik, 
   Robert-Mayer-Str.~10,\\ 
   60054 Frankfurt, 
   Germany
  }
  \and
  \authname{Jan Metzger}
  \authaddress
  {
   Universit\"at Potsdam,
   Institut f\"ur Mathematik,
   Am Neuen Palais 10,
   14469 Potsdam,
   Germany 
  }
}
\date{}

\begin{document}
\hyphenation{}
\pagestyle{footnumber}
\maketitle
\thispagestyle{footnumber}
\begin{abst}%
  We show the existence of a smooth spherical surface minimizing the
  Willmore functional subject to an area constraint in a compact
  Riemannian three-manifold, provided the area is small
  enough. Moreover, we classify complete surfaces of Willmore type
  with positive mean curvature in Riemannian three-manifolds.
\end{abst}
\section{Introduction}
\label{sec:introduction}
For a three-dimensional complete Riemannian manifold $(M,g)$ and an
immersion $f:\Sigma \to M$ the Willmore functional is defined by
\begin{equation*}
  \CW(f) =\frac12 \int_\Sigma H^2 \dmu,
\end{equation*}
where $H$ is the induced mean curvature and $\mu$ the induced area
measure. In the following we let $A$ resp.\ $\Acirc$ be the second
fundamental form resp.\ the trace-free second fundamental form of the
immersion $f$.  Critical points of $\CW$ are called Willmore surfaces
and they are solutions of the Euler-Lagrange equation
\begin{equation*}
  \Delta H +H  |\Acirc|^2 + H \Ric(\nu,\nu)  = 0,
\end{equation*}
where $\Ric$ denotes the Ricci curvature of $(M,g)$ and $\nu$ is the
normal vector to $\Sigma$ in $M$.

In the literature other possible definitions of the Willmore
functional for immersions in a Riemannian manifold were
considered, for example:
\begin{equation*}
  \int_\Sigma |A|^2 \dmu,
  \ \ \
  \int_\Sigma |\Acirc|^2 \dmu,
  \ \ \ \text{or} \ \ \
  \int_\Sigma (H^2+\kappa^M)\dmu.
\end{equation*}
Here $\kappa^M$ denotes the sectional curvature of $M$. In a curved
ambient manifold the Gauss equation and the Gauss-Bonnet Theorem yield
\begin{equation*}
  \CW(\Sigma)
  =
  \frac12 \int_\Sigma |A|^2 \dmu+\int_\Sigma G(\nu,\nu) \dmu+2\pi (1-q(\Sigma))
\end{equation*}
where $q(\Sigma)$ is the genus of $\Sigma$, $G=\Ric-\frac12 \Scal g$
is the Einstein tensor, $\Ric$ denotes the Ricci curvature and $\Scal$
the scalar curvature of $(M,g)$. Hence, the functionals above differ
by lower order terms involving the curvature of $(M,g)$. In
particular, if $(M,g)=(\IR^3,\delta)$ the only difference is a
multiple of $q(\Sigma)$.

In this paper we are interested in surfaces of Willmore type, i.e.\
critical points of $\CW$ subject to an area constraint. These surfaces
are solutions of the Euler-Lagrange equation
\begin{equation*}
  \Delta H + H |\Acirc|^2 + H \Ric(\nu,\nu)  + \lambda H= 0,
\end{equation*}
with the Lagrange parameter $\lambda$.
In \cite{lamm-metzger:2010} we studied spherical surfaces of Willmore
type with positive mean curvature in small geodesic balls. Assuming a
certain lower bound on the Lagrange parameter, we showed that such
surfaces can only concentrate at critical points of the scalar
curvature of $M$. This paper establishes the existence of minimizers
of $\CW$ with fixed small area and verifies that the minimizers
satisfy the assumptions of \cite{lamm-metzger:2010}. Our main result
is as follows:
\begin{theorem}
  \label{main}
  Let $(M,g)$ be a compact, closed Riemannian manifold. Then there
  exists a constant $a_0>0$ such that for all $a\in(0,a_0)$ there is a
  smooth spherical surface $\Sigma_a$ of positive mean curvature that
  minimizes the Willmore functional among all immersed surfaces with
  area $a$.
  
  For any sequence $a_i \to 0$ there is a subsequence $a_{i'}$ such
  that $\Sigma_{a_{i'}}$ is asymptotic to a geodesic sphere centered
  at a point $p\in M$ where $\Scal$ attains its maximum.
\end{theorem}
The existence of a $W^{2,2} \cap W^{1,\infty}$ conformal immersion
minimizing $\CW$ with prescribed small area was recently and
independently obtained by Chen and Li~\cite{CL:2011}.

Kuwert and Sch\"atzle \cite{kuwert-schatzle:2010} constructed smooth
minimizers in a conformal class of the Willmore functional in
$\IR^3$. Existence was recently generalized to arbitrary co-dimension
by Kuwert and Li \cite{kuwert-li:2010} in the class $W^{2,2} \cap
W^{1,\infty}$ and by Rivi{\`e}re \cite{riviere:2010} who also showed
smoothness using his previous results from \cite{riviere:2008}.

The existence of a smooth minimizer of the Willmore functional in
$\IR^n$ with prescribed genus was first proved by Simon
\cite{simon:1993existencewillmore} under a Douglas-type condition,
which was established by Bauer and Kuwert
\cite{bauer-kuwert:2003} (see also \cite{riviere:2010} for a new proof of this result). Recently, Schygulla \cite{Schygulla:2011}
suitably modified the arguments of Simon in order to prove the
existence of a minimizing Willmore sphere in $\IR^3$ with prescribed
isoperimetric ratio.

Under suitable curvature assumptions, Kuwert, Mondino, and Schygulla
\cite{KMS:2011} recently showed the existence of smooth spherical
minimizers of the functionals $\int_\Sigma |A|^2 \dmu$ and
$\int_\Sigma (1+|H|^2)\dmu$ in Riemannian manifolds. The role of the
curvature assumptions is to ensure a uniform bound on the area of
surfaces in a minimizing sequence.

Some other existence and non-existence results of critical points of
$\CW$ by Mondino can be found in \cite{mondino:2010,mondino:2011}.

Previously, in a joint work with F.~Schulze
\cite{lamm-metzger-schulze:2009}, we proved the existence of a
foliation of the end of an asymptotically flat manifold with positive
mass by spherical surfaces of Willmore type. For the existence result
we studied perturbed geodesic spheres and we used the implicit
function theorem to find suitable deformations of these spheres.

In this paper we use the direct method of the calculus of variations
in order to construct the minimizer of $\CW$. A major difficulty is
the invariance of $\CW$ under diffeomorphisms. We overcome this
problem by showing that spherical surfaces with small enough area have
small diameter since the Willmore energy is a priori close to
$8\pi$. Hence the surface is contained in a small geodesic ball around
a point $p\in M$ where the metric $g$ is a small perturbation of the
Euclidean metric. Thus we can apply a result of De Lellis and M\"uller
\cite{delellis-muller:2005,delellis-muller:2006} which gives the
existence of a $W^{2,2}\cap W^{1,\infty}$ conformal parametrization
$F:S^2\to \Sigma \subset \IR^3$ of $\Sigma$.  Therefore, instead of
studying minimizing sequences of immersions $f_k:\Sigma_k\to M$ of
$\CW$, we consider minimizing sequences of parametrizations $F_k\in
W^{2,2}\cap W^{1,\infty}(S^2 ,\IR^3)$. Within this class we are able
to show the existence of a minimizer of $\CW$.

In order to show the higher regularity of the minimizer in
$W^{2,2}\cap W^{1,\infty}$ we suitably modify the arguments of
\cite{KMS:2011} and \cite{Schygulla:2011}.  In our situation their
arguments heavily simplify since the smallness of the area rules out
bad points for the minimizing sequence and we can use the fact that
there exists a limiting parametrization.

As a last step we show that the minimizers we construct satisfy the
assumptions of the main result in \cite{lamm-metzger:2010} and hence
we conclude that the minimizing surfaces have to concentrate around a
maximal point of the scalar curvature of $M$ as the area tends to
zero.

In the following we give a brief outline of the paper. In
section~\ref{sec:preliminaries} we review manifolds with bounded
geometry and show that the diameter estimates of
\cite{simon:1993existencewillmore} extend to these ambient
manifolds. Moreover we show that by a scaling argument one can adjust
the area of a surface without changing the Willmore energy too
much. This fact will be crucial in the proof of the smoothness of the
minimizer of $\CW$ with prescribed area.

In section~\ref{sec:analysis} we study the Willmore functional for
immersions $F\in W^{2,2}\cap W^{1,\infty}(S^2,\IR^3)$ and we show the
lower semi-continuity and the differentiability of the functional.

In section~\ref{sec:direct-minimization} we construct the smooth
minimizer and hence prove the existence part of Theorem~\ref{main}
using the methods described above.

In section~\ref{sec:geom-crit-points} we prove an integral estimate
and and estimate for the Lagrange parameter of the minimizer. These
estimates then allow us to apply the results of
\cite{lamm-metzger:2010}.

In the appendix we derive a variant of the stability inequality of
minimal surfaces for surfaces of Willmore type with positive mean
curvature. Using the methods of \cite{fischer-schoen:1980} we are then
able to classify complete surfaces of Willmore type with positive mean
curvature in Riemannian manifolds. As a corollary we obtain that the
only complete surfaces of this type in $\IR^3$ are round spheres.


%
\section{Preliminaries}
\label{sec:preliminaries}
\subsection{Manifolds with bounded geometry}
\label{sec:unif-norm-coord}
In this section we recall some basic properties of manifolds with
bounded geometry. The main point is that such manifolds have uniformly
controlled normal coordinates.
\begin{definition}
  \label{thm:def_bdd_geometry}
  Let $(M,g)$ be a complete Riemannian manifold. We say that $(M,g)$
  has bounded geometry, if there exists a constant $0<C_B<\infty$ such
  that for each $p\in M$ we have $\inj(M,g,p)\geq C_B^{-1}$ and
  if the Riemann tensor and its first derivative are bounded $|\Riem|+|\nabla \Riem| \leq C_B$.
\end{definition}
In the following we use $B_r$ to denote a Euclidean ball centered at
the origin of radius $r$ and $\CB_r(p)\subset M$ to denote a geodesic
ball of radius $r$ centered at $p\in M$.
\begin{remark}
  \label{rem:uniform_normal_coordinates}
  If $(M,g)$ has bounded geometry with constant $C_B$, then there
  exist constants $h_0<\infty$ and $\rho_0>0$, depending only on
  $C_B$, such that for every $p\in M$, we can introduce normal
  coordinates $\phi:B_{\rho_0} \to \CB_{\rho_0}(p)$ for the metric $g$
  such that in these coordinates the metric $g$ satisfies
  \begin{equation*}
    g = g^E + h, 
  \end{equation*}
  with 
  \begin{equation*}
    \sup_{B_{\rho_0}} (|x|^{-2} |h| + |x|^{-1} |\del h| + |\del^2 h|) \leq h_0.
  \end{equation*}
  Here $g^E$ denotes the Euclidean metric induced by the normal
  coordinates, $|x|$ denotes the Euclidean distance to $p$
  and $\del$ the connection of $g^E$.

  If $(M,g)$ is compact or asymptotically flat, then it is of bounded
  geometry for some constant $C_B$.
\end{remark}

\subsection{Area, diameter and the Willmore energy}
\label{sec:diameter-estimate}
To proceed, we need a generalization of Lemma 1.1 from
\cite{simon:1993existencewillmore} to general ambient manifolds. Although this is
straight forward, we present the proof here to show where the non-flat
ambient geometry has to be taken into account. For ease of
presentation, we split this into three separate statements.
\begin{lemma}[\protect{cf. \cite[Lemma 2.2]{lamm-metzger:2010}}]
  \label{thm:area_leq_diameter}
  Let $g = g^E + h$ on $B_{\rho_0}$ be given, such that
  \begin{equation*}
    \sup_{B_{\rho_0}} (|x|^{-2}|h| + |x|^{-1} |\del h|)
    \leq
    h_0.
  \end{equation*}
  Then there exists a purely numerical constant $c>0$ such that if
  $\rho_1 := \min\{\rho_0, \frac{1}{c\sqrt{h_0}}\}$, then for all surfaces
  $\Sigma \subset B_r$ with $r \in(0,\rho_1)$ we have that
  \begin{equation*}
    |\Sigma|
    \leq
    r^2 \CW(\Sigma).
  \end{equation*}
\end{lemma}
The following proposition is very similar to the calculations in
\cite[Section 1]{simon:1993existencewillmore}. We add only minor modifications in order to deal
with the non-flat background.
\begin{proposition}
  \label{thm:prop_diam_leq_area}
  Let $g = g^E + h$ on $B_{\rho_0}$ be given, such that
  \begin{equation*}
    \sup_{B_{\rho_0}} (|x|^{-2}|h| + |x|^{-1} |\del h|)
    \leq
    h_0.
  \end{equation*}  
  Then there exists a purely numerical constant $c$, such that for
  every smooth surface $\Sigma\subset B_{\rho_0}$ with $\del\Sigma \subset
  \del B_{\rho_0}$ and $0\in\Sigma$ we have that
  \begin{equation*}
    \pi
    \leq
    c \big( (1 + h_0 r^2)  r^{-2} |\Sigma_r| +  \CW(\Sigma_r) \big)
  \end{equation*}
  for all $r\leq \rho$. Here $\Sigma_r := \Sigma\cap B_r$.
\end{proposition}
\begin{proof}
  In $B_{\rho_0}$ consider the position vector field $x$. We denote by
  $c$ a constant that is purely numerical, but which may change from
  line to line. For all surfaces $\Sigma$ as in the statement we have
  \begin{equation*}
    |\div_\Sigma x  - 2| \leq c h_0|x|^2
  \end{equation*}
  in view of the assumption on $h$. Furthermore, we calculate that in
  $B_{\rho_0}$
  \begin{equation*}
    \rmd |x| =  \frac{x}{|x|}.
  \end{equation*}
  In particular, away from the origin, we have
  \begin{equation*}
    \div_\Sigma (|x|^{-2} x)
    =
    |x|^{-2} \div_\Sigma x - 2 |x|^{-3} \rmd |x|(x^T).
  \end{equation*}
  Thus
  \begin{equation}
    \label{eq:1}
    |\div_\Sigma (|x|^{-2} x) - 2|x|^{-4} |x^\perp|^2 | \leq c h_0.
  \end{equation}
  where $x^\perp$ denotes the projection of $x$ onto the normal bundle
  of $\Sigma$.

  Choose $0<s<r<\rho$ such that $\Sigma$ intersects $\del B_r$ and
  $\del B_s$ transversely (note that the set of radii satisfying this
  condition is dense in $(0,\rho)$). Let
  \begin{equation*}
    \Sigma_{s,r} := \Sigma \cap (B_r\setminus \bar B_s)
  \end{equation*}
  and integrate equation~\eqref{eq:1} on $\Sigma_{s,r}$ to obtain
  \begin{equation}
    \label{eq:2}
    \left|\int_{\Sigma_{s,r}} \div_\Sigma (|x|^{-2} x) \dmu
    -
    \int_{\Sigma_{s,r}} 2 |x|^{-4} |x^\perp|^2 \dmu
    \right|
    \leq c h_0 |\Sigma_{s,r}|.    
  \end{equation}
  Using Stokes, we infer that
  \begin{equation}
    \label{eq:3}
    \begin{split}
      &\int_{\Sigma_{s,r}} \div_\Sigma (|x|^{-2} x) \dmu
      \\
      &\quad
      =
      \int_{\Sigma_{s,r}} H |x|^{-2} \la x,\nu\ra\dmu
      -
      \int_{\Sigma\cap \del B_s} |x|^{-2} \la x,\eta\ra \dsig
      +
      \int_{\Sigma\cap \del B_r} |x|^{-2} \la x,\eta\ra \dsig.
    \end{split}
  \end{equation}
  Here $\nu$ is the normal vector of $\Sigma$ and $\eta$ denotes the co-normal of $\del B_s\cap\Sigma$ and $\del
  B_r \cap\Sigma$ in $\Sigma$ respectively. We chose the orientation
  so that $\eta$ points in direction of $\nabla r$. To proceed, we
  note that
  \begin{equation*}
      \int_{\Sigma\cap \del B_\sigma} |x|^{-2} \la x,\eta\ra \dsig
      =
      \sigma^{-2} \int_{\Sigma\cap \del B_\sigma}  \la x,\eta\ra \dsig
      =
      \sigma^{-2} \int_{\Sigma_\sigma} \div_\Sigma x^T \dsig,
  \end{equation*}
  so that 
  \begin{equation}
    \label{eq:4}
    \left| \int_{\Sigma\cap \del B_\sigma} |x|^{-2} \la x,\eta\ra
      \dsig
      -
      2 \sigma^{-2} |\Sigma_\sigma|
      +
      \sigma^{-2} \int_{\Sigma_\sigma} H \la x^\perp,\nu\ra \dmu
    \right|
    \leq c h_0 |\Sigma_\sigma|.
  \end{equation}
  Furthermore
  \begin{equation}
    \label{eq:5}
    \int_{\Sigma_{s,r}}
    |x|^{-4} |x^\perp|^2 - \half H |x|^{-2}\la x,\nu\ra\dmu
    =
    \int_{\Sigma_{s,r}} \big| |x|^{-2}x^\perp - \tfrac14 H\nu\big|^2 -
    \tfrac{1}{16} H^2 \dmu.
  \end{equation}
  Inserting \eqref{eq:3}--\eqref{eq:5} into~\eqref{eq:2}, we infer
  that
  \begin{equation*}
    \begin{split}
      &\int_{\Sigma_{s,r}} \big| |x|^{-2}x^\perp - \tfrac14 H\nu\big|^2 \dmu
      -\half s^{-2} \int_{\Sigma_s} H\la x,\nu\ra \dmu
      + s^{-2} |\Sigma_s|
      \\
      &\quad
      \leq
      r^{-2}|\Sigma_r|
      + \tfrac18\CW(\Sigma_{s,r})
      - \half r^{-2} \int_{\Sigma_r} H\la x,\nu\ra \dmu
      + c h_0 |\Sigma_r|.
    \end{split}
  \end{equation*}
  Since $\Sigma$ is smooth at the origin, we can let $s\to 0$ and drop
  the square term on the left to obtain
  \begin{equation*}
    \pi
    \leq
    r^{-2} |\Sigma_r|
    +
    \tfrac18\CW(\Sigma_r)
    - 
    \half r^{-2} \int_{\Sigma_r} H\la x,\nu\ra \dmu
    + c h_0 |\Sigma_r|.
  \end{equation*}
  Using Cauchy-Schwarz on the third term and recalling that $r<\rho_0$
  we infer the estimate
  \begin{equation}
    \label{eq:6}
    \pi\leq c \big( (1 + h_0 r^2)  r^{-2} |\Sigma_r| +  \CW(\Sigma_r) \big)
  \end{equation}
  for a purely numerical constant $c$. Since the values of $r$ for
  which~\eqref{eq:6} holds are dense in $(0,\rho_0]$ we arrive at the
  claimed estimate by approximation.
\end{proof}
The next lemma shows that the diameter of a surface $\Sigma$ contained
in a Riemannian manifold $(M,g)$ is bounded in terms of its area and
its Willmore energy. We define
\begin{equation*}
  \diam(\Sigma) := \max \{ d_{(M,g)}(p,q) : p,q\in\Sigma\}
\end{equation*}
to be the extrinsic diameter of $\Sigma$. Here $d_{(M,g)}(p,q)$
denotes the geodesic distance of $p$ and $q$ in the ambient manifold
$M$.
\begin{lemma}
  \label{thm:diam_leq_area}
  Let $(M,g)$ be a manifold with $C_B$-bounded geometry. Then there
  exists a constant $C$ depending only on $C_B$ such that for all
  smooth connected surfaces $\Sigma$ we have
  \begin{equation*}
    \diam(\Sigma)
    \leq
    C \big(|\Sigma|^{1/2}\CW(\Sigma)^{1/2} + |\Sigma|\big)
  \end{equation*}
\end{lemma}
\begin{proof}
  Let $h_0$ and $\rho_0$ be as in
  remark~\ref{rem:uniform_normal_coordinates}. Choose $p,q\in\Sigma$
  such that $d:= d_{(M,g)}(p,q) = \diam(\Sigma)$. Assume for now that
  $r\in (0,\tfrac d2)$ is chosen such that $r<\rho_0$.

  Let $N$ be the largest integer smaller than $d/r$ and let $p_0=p$. For
  $j=1,\ldots,N-1$ we choose $p_j\in\Sigma$ at distance $(j+\half)r$
  to $p_0$, which is possible since $\Sigma$ is connected. Then the geodesic balls $\CB_{r/2}(p_j)$ are pairwise disjoint
  for $j=0,\ldots,N-1$. Using proposition~\ref{thm:prop_diam_leq_area}
  with $p_j$ as center and summing over $j$ yields
  \begin{equation}
    \label{eq:7}
    N\pi \leq c\big(\CW(\Sigma) + (1 + h_0 r^2) r^{-2}|\Sigma|\big).    
  \end{equation}
  With this in mind, we let 
  \begin{equation*}
    r := \min \left\{ \rho_1/2, \frac{1}{4}
      \sqrt{\frac{|\Sigma|}{\CW(\Sigma)}} \right\},
  \end{equation*}
  where $\rho_1 = \min\{\rho_0, \frac{1}{c\sqrt{h_0}} \}$ is such that
  lemma~\ref{thm:area_leq_diameter} applies with $\rho_0$ and $h_0$ as
  above.

  We have to check that $r<d/2$. Assume for the contrary that $d/2
  \leq r$. Then in particular $d\leq \rho_1$ and
  lemma~\ref{thm:area_leq_diameter} implies that
  \begin{equation*}
    \sqrt{\frac{|\Sigma|}{\CW(\Sigma)}}
    \leq
    d
    \leq
    2 r
    \leq
    \frac12 \sqrt{\frac{|\Sigma|}{\CW(\Sigma)}},
  \end{equation*}
  a contradiction. Hence $r< d/2$ and thus $N\geq
  \tfrac{d}{2r}$. Revisiting~\eqref{eq:7} thus yields
  \begin{equation*}
    d
    \leq
    c \big( r\CW(\Sigma) +  (1+h_0 r^2) r^{-1} |\Sigma| \big)
  \end{equation*}
  and since
  \begin{equation*}
    r^{-1} \leq \max \left\{ \frac{2}{\rho_1}, 4
      \sqrt{\frac{\CW(\Sigma)}{|\Sigma|}}\right\}
  \end{equation*}
  we find that
  \begin{equation*}
    \begin{split}
      d
      &\leq
      c (1+ r^2 h_0) \sqrt{|\Sigma|\CW(\Sigma)} + c/\rho_1 |\Sigma|
      \\[1ex]
      &\leq
      c (1+ r^2 h_0)  \sqrt{|\Sigma|\CW(\Sigma)} + c(\rho_0^{-1} + \sqrt{h_0}) |\Sigma|.
    \end{split}
  \end{equation*}
  This yields the claimed estimate.
\end{proof}
\subsection{Area adjustment by scaling}
\begin{lemma}
  \label{thm:area_adjust}
  Let $(M,g)$ be a manifold with $C_B$-bounded geometry. Then there
  exists $\rho_1>0$ with the following property. Let $r\in(0,\rho_1)$, 
  $p\in M$, and let $x$ be the position vector field with respect to
  geodesic normal coordinates in $\CB_r(p)$. Denote by $\Phi :
  \CB_{r/4}(p) \times (-\infty,2)  \to \CB_r(p)$ the flow associated to
  $x$.
  Then for every $a\in \IR$ and $\Sigma\subset \CB_{r/4}(p)$ with
  $|\Sigma| \in (\frac{a}{2},\frac{3a}{2})$ there exists $t_0\in \IR$
  with $|\Phi_{t_0}(\Sigma)| = a$ and $|t_0| \leq 2 \frac{| |\Sigma| - a |}{a}$.
\end{lemma}
\begin{proof}
  Let $\rho_0$ and $h_0$ as in
  Remark~\ref{rem:uniform_normal_coordinates}. We choose
  $\rho_1\in(0,\rho_0)$ such that
  \begin{equation*}
    |\nabla x - \Id | \leq \tfrac{1}{2}
  \end{equation*}
  on $\CB_{\rho_1}(p)$ for all $p\in M$. Note that then $\rho_1$
  depends only on $C_B$.
  
  Let $r\in(0,\rho_1)$ and $\Sigma'\subset \CB_{r}(p)$ be an arbitrary
  surface with $|\Sigma'|\geq a/2$. Then we calculate that
  \begin{equation*}
    \DDeval{}{t}{t=0} |\Sigma'|
    =
    \int_{\Sigma'} H \la x, \nu\ra \dmu
    =
    \int_{\Sigma'} \div_{\Sigma'} x \dmu
    \geq
    |\Sigma'|
    \geq
    \frac{a}{2}.
  \end{equation*}
  If we consider $\Sigma$ as in the statement of the lemma, we can
  apply this estimate to $\Sigma_t= \Phi_t(\Sigma)$ as long as
  $\Sigma_t \in B_r(p)$ and $|\Sigma_t|\in
  (\frac{a}{2},\frac{3a}{2})$. In particular the area of $\Sigma_t$ is
  a continuous and strictly increasing function of $t$. In addition we
  have that
  \begin{equation*}
    \begin{split}
      |\Sigma_{t^+}| &\geq a
      \qquad\text{for}\qquad
      t^+ := \max\big\{ 0, 2 \tfrac{a-|\Sigma|}{a} \big\}, \text{ and}
      \\
      |\Sigma_{t^-}| &\leq a
      \qquad\text{for}\qquad
      t^- := \min\big\{ 0, 2 \tfrac{|\Sigma|-a}{a} \big\}.
    \end{split}
  \end{equation*}
  This yields the claim.
\end{proof}
\begin{lemma}
  \label{thm:curvature_scaling}
  Let $(M,g)$ be a manifold with $C_B$-bounded geometry and let $\rho_1$,
  be as in Lemma~\ref{thm:area_adjust}. There exists a constant $C$
  depending only on $C_B$ with the following property. Let
  $r\in(0,\rho_1)$ and let $\Sigma\subset \CB_{r/2}(p)$. Then
  \begin{equation*}
    \DDeval{}{t}{t=0} \int_{\Sigma_t} |A_t|^2 \dmu_t
    \leq    
    C |\Sigma|^{1/2} \left(\int_\Sigma H^2\dmu\right)^{1/2}
    + Cr \int_\Sigma |A|^2 \dmu
    + C(1+r) |\Sigma|. 
  \end{equation*}
  Here we use the notation of Lemma~\ref{thm:area_adjust}, so that
  $\Sigma_t=\Phi_t(\Sigma)$, $A_t$ denotes the second fundamental form
  of $\Sigma_t$ and $\dmu_t$ its induced measure.
\end{lemma}
\begin{proof}
  We use the Gauss equation to write 
  \begin{equation*}
    \int_\Sigma |A|^2 \dmu
    =
    \int_\Sigma H^2 \dmu - 2\int_\Sigma G(\nu,\nu)\dmu - 4\pi(1-g(\Sigma))
  \end{equation*}
where $g(\Sigma)$ is the genus of $\Sigma$ and $G(\nu,\nu)=\Ric-\frac12 \Scal g$. Therefore we have
  \begin{equation*}
    \DDeval{}{t}{t=0} \int_{\Sigma_t} |A_t|^2 \dmu_t
    =
   2  \delta_{\la x,\nu\ra} \CW(\Sigma) - \delta_{\la x,\nu\ra} \CV(\Sigma)     
  \end{equation*}
  where
  \begin{equation*}
    \CV(\Sigma) = 2 \int_\Sigma G(\nu,\nu) \dmu.
  \end{equation*}
  We estimate the variations of $\CW$ and $\CV$ separately.  We have
  that
  \begin{equation}
    \label{eq:33}
    | \nabla x - \Id | \leq C |x|^2
    \qquad\text{and}\qquad
    |\nabla^2 x | \leq C
  \end{equation}
  where $C$ is a constant depending only on $C_B$.

  Consider the function $\la x,\nu\ra$ on $\Sigma$. A calculation
  shows that with respect to an adapted ON-Frame $\{e_1,e_2,e_3
  =\nu\}$, we have
  \begin{equation*}
    \begin{split}
      \Delta \la x,\nu \ra
      &=
      \la \nabla_{e_i,e_i} x, \nu \ra
      - H \la \nabla_\nu x,\nu\ra
      + 2 \la \nabla_{e_i} x, e_k\ra A_{ik}\\
      &\phantom{=}
      - \la x,\nu\ra |A|^2
      + \la x,  e_k\ra \nabla_i A_{ik}
      \\
      &=
      H
      -\la x,\nu\ra |A|^2
      + \la x,e_k\ra \nabla_{e_k}H
      + O(1) + O(r) * A.
    \end{split}
  \end{equation*}
  In the last equality we used the Codazzi equation to rewrite $\div A
  = \nabla H + \Ric(\nu,\cdot)$ together with the fact that $x=
  O(r)$. In addition we used~\eqref{eq:33} and use the notation $O(1)$
  for terms which are bounded by a constant $C$ and $O(r) * A$ for
  terms bounded by $Cr|A|$. Here as ususal $C$ depends only on $C_B$.

  We calculate the variation of the Willmore functional with respect
  to scaling:
  \begin{equation}
    \label{eq:26}
    \begin{split}
      \delta_{\la x,\nu\ra} \CW
      &=
      \int_\Sigma \la x,\nu\ra \big( \Delta H + H | \Acirc |^2 + H
      \Ric(\nu,\nu) \big) \dmu
      \\
      &=
      \int_\Sigma H \Delta \la x,\nu \ra + \la x,\nu\ra H| \Acirc|^2 +
      \la x, \nu\ra H \Ric(\nu,\nu) \dmu
      \\
      &=
      \int_\Sigma H^2 - \half H^3 \la x,\nu\ra + H \nabla_k H \la x,
      e_k\ra
      + H \big( O(1) + O(r) * A\big) \dmu.
    \end{split}
  \end{equation}
  Integrate by parts the third term on the right to calculate further
  \begin{equation*}
    \begin{split}
      \int_\Sigma H \nabla_k H \la x, e_k\ra \dmu
      &=
      \int_\Sigma \half \nabla_k (H^2) \la x, e_k\ra\dmu
      =
      - \half \int_\Sigma H^2 \nabla_k \la x, e_k\ra \dmu
      \\
      &=
      - \half \int_\Sigma H^2  \div_\Sigma x - H^3\la x, \nu\ra\dmu.
    \end{split}
  \end{equation*}
  Inserting this into equation~\eqref{eq:26} and taking into
  account~\eqref{eq:33} yields that
  \begin{equation*}
    \delta_{\la x,\nu\ra} \CW = \int_\Sigma H \big( O(1) + O(r) * A\big) \dmu.
  \end{equation*}
  This can be estimated as follows:
  \begin{equation*}
    |\delta_{\la x,\nu\ra} \CW|
    \leq
    C |\Sigma|^{1/2} \left( \int_\Sigma H^2 \dmu \right)^{1/2}
    +
    C r \int_\Sigma |A|^2 \dmu.
  \end{equation*}
  We proceed with the variation of $\CV(\Sigma)$. From
  \cite[(75)]{lamm-metzger-schulze:2009} we obtain that
  \begin{equation*}
    \tfrac{1}{2}\delta_{\la x,\nu\ra} \CV(\Sigma)
    =
    \int_\Sigma \la x, \nu \ra \big( \nabM_\nu G(\nu,\nu) + H
    G(\nu,\nu) \big) - 2 G\big( \nu,\nabla(\la x,\nu\ra) \big) \dmu.
  \end{equation*}
  Straight forward estimates show that
  \begin{equation*}
    |\delta_{\la x,\nu\ra} \CV(\Sigma)|
    \leq
    C (1+r) |\Sigma|
    +
    C r |\Sigma|^{1/2} \left(\int_\Sigma |A|^2 \dmu \right)^{1/2}.
  \end{equation*}
  This implies the claim.
\end{proof}
\begin{lemma}
  \label{thm:area_adjust_with_curvature}
  Let $(M,g)$ be a manifold with $C_B$-bounded geometry and let
  $\rho_1$, be as in Lemma~\ref{thm:area_adjust}. For every constant
  $C_0$ there exists a constant $C_1$ with the following
  properties. If $r\in(0,\rho_1)$, $a\in(0, C_0r^2)$ and
  $\Sigma\subset B_{r/4}(p)$ with $|\Sigma| \in (\frac{a}{2} ,
  \frac{3a}{2})$ then there exists a surface $\Sigma'\subset \CB_{r}(p)$
  with $|\Sigma'| = a$ and
  \begin{equation*}
    \left| \int_{\Sigma'}|A|^2 \dmu - \int_{\Sigma}|A|^2 \dmu \right|
    \leq
    C_1 r \frac{| |\Sigma| - a |}{a}
    \left( 1 + \int_{\Sigma}|A|^2 \dmu \right).    
  \end{equation*}
\end{lemma}
\begin{proof}
  Let $\Sigma$ be as in the statement. Using
  Lemma~\ref{thm:area_adjust} we can find $\Sigma'\subset \CB_r(p)$ with
  $|\Sigma'|=a$ in the form $\Sigma' = \Phi_{t_0}(\Sigma)$, where
  $\Phi_t$ is as in Lemma~\ref{thm:area_adjust}. In addition, we have
  that $|t_0| \leq 2 \frac{| |\Sigma| -a |}{a} \leq 2$.
  
  To analyze the amount the second fundamental form has changed, we
  assume for definiteness that $t_0>0$. From
  Lemma~\ref{thm:curvature_scaling} we find that for all $\Sigma_t =
  \Phi_t(\Sigma)$ with $t\in[0,t_0]$ we have that
  \begin{equation*}
    \DD{}{t} \left(1 +\int_{\Sigma_t} |A_t|^2 \dmu_t\right)
    \leq
    C r \left(1 + \int_{\Sigma_t} |A_t|^2 \dmu_t \right),
  \end{equation*}
  where the constant $C$ only depends on $C_B$, $C_0$ and $\rho_1$. In
  particular, we used that the area of all $\Sigma_t$ is bounded by
  $\frac{3C_0}{2} r^2 \leq \frac{3C_0}{2} \rho_1 r$. Integrating this
  ordinary differential inequality on $[0,t_0]$ and using the fact
  that $|t_0|\leq 2$ is a priori bounded, we arrive at the claimed
  estimate.
\end{proof}


%
\section{Analytical aspects of the Willmore functional}
\label{sec:analysis}
In this section we consider the Willmore functional in the space of
pa\-ra\-met\-ri\-za\-tions which are in a subset of
$W^{2,2}(S^2,\IR^3) \cap W^{1,\infty}(S^2,\IR^3)$. We assume that
$\IR^3$ is equipped with a smooth metric $g$ of which we assume that
with respect to standard coordinates all components and derivatives up
to second order thereof are bounded. When we refer to coordinates on
$\IR^3$ we use Greek indices, and when referring to coordinates on
$S^2$ we use Latin indices.

We define the space
\begin{equation*}
  B
  :=
  \big\{ F \in W^{2,2}(S^2,\IR^3) \cap W^{1,\infty}(S^2,\IR^3)
  \mid
  \bar g,\bar g^{\#} \in L^{\infty}(S^2)\cap W^{1,2}(S^2) \big\}
\end{equation*}
where we denote by $\bar g$ the pull-back metric $F^*g$ on $S^2$ and
by $\bar g^{\#}$ its inverse. The function spaces and all tensor norms
are defined with respect to the standard metric on the sphere and the
standard metric $g^E$ on $\IR^3$.

\subsection{Definition of the Willmore functional}
First we establish that the Willmore functional
\begin{equation*}
  \CW(F) = \frac12 \int_{F(S^2)} H^2 \dmu 
\end{equation*}
is well defined for $F\in B$. Denote by $h_{ij}$ the second
fundamental form of $F(S^2)$, by $\bar\Gamma$ the Christoffel symbols
of $\bar g$ and by $\Gamma$ the Christoffel symbols of the ambient
metric $g$. The Weingarten equation gives
\begin{equation}
  \label{eq:23}
  F_{ij}^\al
  =
  -h_{ij}\nu^\al + \bar\Gamma_{ij}^m F_m^\al - \Gamma_{\be\ga}^\al
  F_i^\be F_j^\ga.
\end{equation}
Here we use the shorthand notation 
\begin{equation*}
  F_i^\al = \dd{F^\al}{x^i}
  \quad\text{and}\quad
  F_{ij}^\al = \dd{^2 F^\al}{x^i \del x^j}.
\end{equation*}
In fact we can take equation~\eqref{eq:23} as the definition of the
second fundamental form. Note that since $\bar g,\bar g^\# \in
W^{1,2}(S^2)$ we have that $\bar\Gamma \in L^2(S^2)$ so that the
second fundamental form is in $L^2(S^2)$. Taking the trace
of~\eqref{eq:23} gives
\begin{equation}
 \label{eq:31}
  H = -(\bar g^{ij} F_{ij}^\al + (\Gamma_{\be\ga}^\al\circ F) F_i^\be
  F_j^\ga)(g_{\al\de}\circ F) \nu^\de.
\end{equation}
In particular, we can write the Willmore functional as
\begin{equation}
  \label{eq:32}
  \CW(F)
  =
  \frac12 \int_{S^2}
  \big[(\bar g^{ij} F_{ij}^\al + (\Gamma_{\be\ga}^\al\circ F) F_i^\be
  F_j^\ga)(g_{\al\de}\circ F) \nu^\de\big]^2 \sqrt{|\bar g|} \dx  
\end{equation}
where $|\bar g| = \det(\bar g)$ and $\dx$ denotes the standard volume
element on $S^2$. Clearly $\CW$ is continuous on $B$ where we equip
$B$ with the topology induced by convergence of $F$ in
$W^{2,2}(S^2,\IR^3) \cap W^{1,\infty} (S^2,\IR^3)$ and $\bar g$ and
$\bar g^\#$ in $W^{1,2}\cap L^{\infty}$.

\subsection{Lower semi-continuity}
In this section we show lower semi-continuity of $\CW$ in $B$ with
respect to weak convergence.
\begin{proposition}
  \label{thm:lower-semi-continuity}
  Assume that $F_k,F \in B$ are parametrizations of surfaces
  $\Sigma_k$ and $\Sigma$ and that $g$ is a smooth metric on
  $\IR^3$ such that all coefficients of $g$ with respect to standard
  coordinates on $\IR^3$ and all their derivatives are bounded.
  If 
  \begin{equation}
    \label{eq:29}
    \begin{cases}
      F_k \weakto F &\text{weakly in $W^{2,2}(S^2,\IR^3,g^E)$\quad and } \\
      F_k \weaktostar F &\text{weakly-* in $W^{1,\infty}(S^2,\IR^3,g^E)$,}
    \end{cases}
  \end{equation}
  then
  \begin{equation*}
    \CW(F) \leq \liminf_{k\to\infty} \CW(F_k).
  \end{equation*}
\end{proposition}
\begin{proof}
  Let $(\bar g_k)_{ij} : = g(\dd{F_k}{x^i},\dd{F_k}{x_j})$ be the
  coefficients of the induced metric on $\Sigma_k$, pulled back to
  $S^2$, $|\bar g_k| = \det (g_k)_{ij}$, and $\nu_k$ the normal of
  $\Sigma_k$ with respect to $g$. The above convergence~\eqref{eq:29}
  implies that for any $1<q<\infty$ we have
  \begin{equation}
    \label{eq:30}
    \begin{aligned}
      (\bar g_k)_{ij} &\to \bar g_{ij}  \text{ in $L^q(S^2)$},\qquad
      &
      (\bar g_k)_{ij} &\weaktostar \bar g_{ij}  \text{ in
        $L^\infty(S^2)$},
      \\
      (\bar g_k)^{ij} &\to \bar g^{ij}  \text{ in $L^q(S^2)$},
      &
      (\bar g_k)^{ij} &\weaktostar \bar g^{ij} \text{ in
        $L^\infty(S^2)$},
      \\
      \nu_k &\to \nu  \text{ in $L^q(S^2)$},
      &
      \nu_k &\weaktostar \nu   \text{ in $L^\infty(S^2)$ },
      \\
      |\bar g_k| &\to |\bar g| \text{ in $L^q(S^2)$},
      &
      |\bar g_k| &\weaktostar |\bar g| \text{ in $L^\infty(S^2)$}.
    \end{aligned}
  \end{equation}
  From equation~\eqref{eq:32} we find that 
  \begin{equation*}
    \CW(F_k)
    =
    \frac12 \int_{S^2}  \left( \bar g_k^{ij} \left( (F_k)^\al_{ij}
        + \Gamma^\alpha_{\beta\gamma}\circ F_k (F_k)^\be_i (F_k)^\ga_j\right) (g\circ
      F_k)_{\alpha\beta} \nu_k^\beta \right)^2 \sqrt{|g_k|} \dx.
  \end{equation*}
  We split this expression into two parts 
  \begin{equation*}
    \CW(F_k)
    =
    \CW_1(F_k) + \CW_2(F_k)
  \end{equation*}
  with
  \begin{equation*}
    \CW_1(F_k)
    =
    \frac12 \int_{S^2}
    \left( \bar g_k^{ij} \dd{^2 F_k^\al}{x^i\del x^j}
      (g\circ F_k)_{\alpha\beta} \nu_k^\beta \right)^2 \sqrt{|g_k|} \dx
  \end{equation*}
  and
  \begin{equation*}
    \begin{split}
      \CW_2(F_k)
      &=
      \frac12 \int_{S^2}
      \Big(
      2 \bar g_k^{ij} (F_k)^\alpha_{ij} (g_{\al\be}\circ
      F_k)\nu_k^\be  \bar g_k^{ab} (\Gamma_{\eps\mu}^\delta \circ F_k)
      (F_k)^\eps_a (F_k)^\mu_b (g_{\delta\rho}\circ F_k) \nu_k^\rho
      \\
      &\qquad
      +
      \left(\bar g_k^{ij} \Gamma_{\be\ga}^\al\circ F_k
        (F_k)^\be_i (F_k)^\ga_j (g_{\al\be}\circ F_k)
        \nu_k^\be\right)^2
      \Big) \sqrt{|g_k|} \dx
    \end{split}
  \end{equation*}
  By the Sobolev embedding and smoothness of $g$ it follows that
  $g_{\al\be}\circ F_k$ converges in $L^\infty(S^2)$ to
  $g_{\al\be}\circ F$ and that $\Gamma^\al_{\be\ga}\circ F_k$
  converges in $L^\infty(S^2)$ to $\Gamma^\al_{\be\ga}\circ F$. In
  view of the convergence~\eqref{eq:30} it thus follows that $\CW_2$
  is continuous with respect to the convergence in~\eqref{eq:29}:
  \begin{equation*}
    \CW_2(F_k) \to \CW_2(F)
    \quad\text{for}\quad
    k\to\infty.
  \end{equation*}
  To analyze $\CW_1(F_k)$ we let $\phi$ be the integrand in the
  definition of $\CW_1(F_k)$:
  \begin{equation*}
    \CW_1(F_k)
    =
    \frac12 \int_{S^2}
    \phi\big(\bar g_k^{ij}, (g_{\al\be}\circ F_k), \nu_k, \sqrt{|g_k|},
    (F_k)^\al_{ij}\big)
    \dx.
  \end{equation*}
  Then $\phi$ is smooth with respect to all variables, $\phi\geq 0$,
  and $\phi$ is convex with respect to the last set of variables
  $(F_k)^\al_{ij}$ as it is the concatenation of the following three
  maps: The linear (in $(F_k)^\al_{ij}$) map
  \begin{equation*}
    \big(\bar g_k^{ij}, (g_{\al\be}\circ F_k), \nu_k, \sqrt{|g_k|},
    (F_k)^\al_{ij}\big)
    \mapsto
    \bar g_k^{ij} (F_k)^\al_{ij} (g_{\alpha\beta}\circ F_k) \nu_k^\beta,
  \end{equation*}
  the convex map $\xi \mapsto \xi^2$ for $\xi\in\IR$ and the linear
  multiplication by $\sqrt{|g_k|}$. Lower semi-continuity then
  follows as in the proof of \cite[Theorem 1.6]{struwe:2008}.

\end{proof}

\subsection{Differentiability}
Given a map $F\in B$ and a smooth vector field $X\in\CX(\IR^3)$ we have that
$X\circ F$  is in $W^{2,2}\cap W^{1,\infty}(S^2,\IR^3,g^E)$. Furthermore,
for small enough $\eps>0$ the map $F_t := F+t(X\circ F)$ is in $B$ for all
$t\in (-\eps,\eps)$ since the inverse of the metric $\bar g_t$ induced
by $F_t$ is a smooth function of $\bar g_t$.

We can thus calculate the variation of $\CW$ in direction of $X$
\begin{equation*}
  \delta_X \CW(F)
  =
  \ddeval{}{t}{t=0} \CW(F_t).
\end{equation*}
A fairly long but standard calculation shows that $\CW$ is indeed
differentiable and that its variation in direction $X$ is given by
\begin{equation*}
  \begin{split}
    \delta_X\CW(F)
    &=
    \int_{S^2}
    - H\bar g^{ij}g(\nabla^2_{F_i,F_j} X,\nu)
    -2 H \bar g^{ik}\bar g^{jl} h_{ij} g(\nabla_{F_i} X, F_j)
    \\ 
    &\phantom{=\int_{S^2}}
    + H^2 g(\nabla_\nu X, \nu)
    - H \Ric(X,\nu)
    + \half H^2 \div^T X \dmu
  \end{split}
\end{equation*}
where $\Ric$ denotes the Ricci curvature of $g$, $\nabla$ the
connection of $g$ and $\nabla^2$ the second covariant derivative of
vector fields with respect to $g$.  There is also a formulation of the
Euler-Lagrange equation for variations $X\in W^{2,2}\cap
W^{1,\infty}(S^2,\IR^3)$ which are not induced by a smooth ambient
variation. To consider such vector fields, derivatives of $X$ have to
be calculated with the pull-back $\nabla^* = F^* \nabla$ of the of the
Levi-Civita connection of $g$. This affects only the way in which the
the second derivatives are calculated. We obtain the following
expression:
\begin{equation*}
  \begin{split}
    \delta_X\CW(F)
    &=
    \int_{S^2}
    - H\bar g^{ij}g((\nabla^*)^2_{F_i,F_j} X,\nu)
    -2 H \bar g^{ik}\bar g^{jl} h_{ij} g(\nabla^*_{F_i} X, F_j)
    \\ 
    &\phantom{=\int_{S^2}}
    - H \Ric(X,\nu)
    + \half H^2 \div^T X \dmu.
  \end{split}
\end{equation*}
Note that if $F$ is $C^4(S^2,\IR^3)$, we can integrate by parts the
terms involving derivaties of $X$ and write the variation in a more
familiar form:
\begin{equation*}
  \delta_X\CW(F)
  =
  \int_{S^2} g(X,\nu) (\Delta H + H |\Acirc|^2 + H \Ric(\nu,\nu)) \dmu.    
\end{equation*}

\section{Direct minimization}
\label{sec:direct-minimization}
In this section we construct minimizers for the Willmore functional
subject to a small area constraint by direct minimization. We assume
that $(M,g)$ is compact without boundary.

Fix a point $p\in M$. For $r<\inj(M,g,p)$ we consider the geodesic
spheres $S_r(p)$. By \cite{mondino:2010} these surfaces satisfy
\begin{equation}
  \label{eq:27}
  \CW(S_r(p))
  =
  8\pi - \frac{4\pi r^2}{3}\Scal(p) + O(r^3).
\end{equation}
In particular, for a given $\eps>0$ there exists a constant $0<a_0=a_0(\eps)$ such
that for $|S_r(p)|\leq a_0$ we have $\CW(S_r) \leq 8\pi + \eps$.

Fix some $a\in(0,a_0)$. We consider a minimizing sequence for $\CW$ of
surfaces $\Sigma_k$ with $|\Sigma_k|=a$. By comparison with geodesic
spheres, we can assume that $\CW(\Sigma_k) \leq 8\pi +\eps$. Thus, in
view of Lemma~\ref{thm:diam_leq_area} there exists a constant $C$ such
that
\begin{equation*}
  \diam(\Sigma_k)
  \leq
  C (a_0^{1/2} + a_0)
\end{equation*}
for all $k$. By choosing $a_0$ small enough, we can ensure that
$4\diam(\Sigma_k) < \inj(M,g)$ uniformly. By compactness of $M$, by
choosing $a_0$ small enough, and by passing to a sub-sequence if
necessary, we can assume that all the $\Sigma_k$ are contained in
$\CB_{\rho_0/16}(p)$ for some suitable $p\in M$, where $\rho_0$ is as
in Remark~\ref{rem:uniform_normal_coordinates}. We decorate all
geometric quantities on $\Sigma_k$ by the sub-script $k$, ie. $H_k,
\nu_k, \ldots$

By the Gauss equation, we have 
\begin{equation*}
  \CW(\Sigma_k)
  =
  8\pi
  +
  \int_{\Sigma_k} |\Acirc_k|^2 \dmu_k
  +
  \int_{\Sigma_k} G(\nu_k,\nu_k) \dmu_k.
\end{equation*}
By assumption we have $\CW(\Sigma_k) \leq 8\pi + \eps$. Since the
curvature of $(M,g)$ is bounded, we can estimate the last term by $Ca_0$.
Thus we obtain the estimate
\begin{equation*}
  \|\Acirc_k\|_{L^2(\Sigma_k)}^2
  \leq
  \eps + Ca_0.
\end{equation*}
From $|A_k|^2 = |\Acirc_k|^2 + \half H_k^2$, we also get
\begin{equation*}
  \|A_k\|_{L^2(\Sigma_k)}^2
  \leq
  8\pi + 2\eps + Ca_0.
\end{equation*}
in the following proposition we show that we can pass this sequence to
a (weak) limit, and that $\CW$ is lower semi-continuous under this
limit.
\begin{proposition}
  Let $(M,g)$ be compact without boundary. Then there exists $\eps>0$,
  depending only on the geometry of $M$, such that the following
  holds.  If $\Sigma_k$ is a sequence of immersed surfaces with
  \begin{equation}
    \label{eq:22}
    |\Sigma_k| = a
    < \eps
    \quad\text{and}\quad
    \CW(\Sigma_k) < 8\pi + \eps
  \end{equation}
  then there exists a family of parametrizations
  \begin{equation*}
    F_k : S^2 \to \Sigma_k
  \end{equation*}
  such that the $F_k$ converge weakly in $W^{2,2}$ and weakly-* in
  $W^{1,\infty}$ to a limiting parametrization
  \begin{equation*}
    F: S^2 \to \Sigma \subset (M,g)
  \end{equation*}
  such that $|\Sigma| =a$ and
  \begin{equation*}
    \CW(\Sigma)
    \leq
    \liminf_{k\to\infty} \CW(\Sigma_k).
  \end{equation*}
\end{proposition}
\begin{proof}
  Let $(\rho_0,h_0)$ be as in
  Remark~\ref{rem:uniform_normal_coordinates}.  Using the reasoning
  prior to the statement of this proposition, we can assume without
  loss of generality, that $\Sigma_k \subset \CB_{\rho_0}(p)$ for some
  $p\in M$. Introducing normal coordinates $x: B_{\rho_0} \to
  \CB_{\rho_0} (p)$ we can pull-back the $\Sigma_k$ for all $k$ and
  the metric $g$ to $B_{\rho_0}$ and $g$ has the form
  \begin{equation*}
    g = g^E + h
  \end{equation*}
  with $h$ as in Remark~\ref{rem:uniform_normal_coordinates}. Assuming
  that $a_0$ is small enough, it is easy to see that
  equation~\eqref{eq:22} implies that $\Sigma_k$ satisfies
  \begin{equation*}
    \CW^E(\Sigma_k) < 8\pi + 2\eps
  \end{equation*}
  where $\CW^E$ denotes the Willmore functional computed with respect to
  the Euclidean background metric $g^E$. Via the Gauss equation with
  respect to the Euclidean background, this implies that
  \begin{equation*}
    \|\Acirc_k^E\|_{L^2(\Sigma_k,g^e)} < 2\eps
  \end{equation*}
  on $\Sigma$. The estimates of DeLellis and Müller
  \cite{delellis-muller:2005,delellis-muller:2006} imply that there
  exist conformal (with respect to the Euclidean background)
  parametrizations $F_k : S^2 \to \Sigma_k$ which are uniformly
  bounded in $W^{2,2}(S^2,\IR^3) \cap W^{1,\infty}(S^2,\IR^3)$. Thus
  there is a subsequence of the $F_k$ which we relabel to $F_k$, such
  that for any given $1\le p <\infty$ we have
  \begin{equation}
    \label{eq:28}
    \begin{cases}
      F_k \weakto F &\text{weakly in $W^{2,2}(S^2,\IR^3,g^E)$} \\
      F_k \weaktostar F &\text{weakly-* in $W^{1,\infty}(S^2,\IR^3,g^E)$,
        \quad and} \\
      F_k \to F & \text{ in $W^{1,p}(S^2,\IR^3,g^E)$}
    \end{cases}
  \end{equation}
  for a function $F\in W^{1,\infty}\cap W^{2,2}(S^2,\IR^3)$. Denote
  $\Sigma = F(S^2)$, then this convergence implies that $|\Sigma| =
  \lim_{k\to\infty} |\Sigma_k| =
  a$. Proposition~\ref{thm:lower-semi-continuity} implies the lower
  semi-continuity of $\CW$ with respect to the convergence in
  equation~\eqref{eq:28}.
\end{proof}

In the following we show that the limiting parametrization $F\in
W^{2,2}\cap W^{1,\infty}(S^2,M)$ of $\Sigma$ of the minimizing sequence
$F_k$ for the Willmore functional is smooth. This follows from
suitable modifications of a recent result of Kuwert-Mondino-Schygulla
\cite{KMS:2011} on the existence of smooth spheres minimizing
$\int_\Sigma |A|^2 \dmu$ resp.\ $\int_\Sigma(|H|^2+1)\dmu$ in
Riemannian manifolds satisfying suitable curvature conditions and a
result of Schygulla \cite{Schygulla:2011} on the existence of a
minimizing Willmore sphere with prescribed isoperimetric ratio in
$\IR^3$. Both of these results rely on the fundamental existence
result (and especially the approximate graphical decomposition lemma)
of Simon \cite{simon:1993existencewillmore}.

In the following we indicate how to modify the arguments of section 3
in \cite{KMS:2011} in order to handle the present situation. First of
all we note that because of the above proposition we get that the
Radon measures on $M$
\begin{equation*}
  \mu_k(E)
  =
  \int_{F_k^{-1}(E)}  \dmu_{S^2}(y)
  \qquad\text{and}\qquad
  \alpha_k(E)= \int_{F_k^{-1}(E)} |A_k|^2 \dmu_k
\end{equation*}
converge weakly to limiting Radon measures $\mu$ and $\alpha$. Note
that for $\CW<8\pi+\eps$ the monotonicity formula implies that the
density of $\mu$ is one on its support. We have that
\begin{equation*}
  \mu(E)=\int_{F^{-1}(E)} \dmu_{S^2}(y)
\end{equation*}
is the induced measure of the limiting immersion $F:S^2 \to \IR^3$.

We define bad points by
\begin{equation*}
  B^\delta =\{ \xi \in  \operatorname{spt} \mu \mid \alpha(\{\xi\})\ge \delta^2\} 
\end{equation*}
and we note that for $\eps$ (and hence $a$) small enough
$B^\delta=\emptyset$. This can be seen as follows:

Assume that $B^\delta\supset \{p_1,\ldots,p_l\}$ and choose a radius $\rho>0$ such that $B_\rho(p_i) \cap B_\rho(p_j)=\emptyset$ for all $1\le i,j \le l$. Then we have
\begin{equation*}  
  \begin{split}
    8\pi+\eps
    &> \lim_{k\to \infty} \CW(\Sigma_k)
    \\
    &\ge \lim_{k\to \infty} \sum_{i=1}^l \CW (\Sigma_k,B_\rho(p_i))
    +\CW(\Sigma_k,\Sigma_k \backslash \cup_{i=1}^l B_\rho(p_i))
    \\
    &\ge
    \CW(\Sigma,\Sigma \backslash \cup_{i=1}^l B_\rho(p_i)) +l\delta^2,
  \end{split}
\end{equation*}
where $\CW(\Sigma,E)=\frac12 \int_{\Sigma \cap E} |H|^2 \dmu$. Since
$H\in L^2(\Sigma)$ we can choose $\rho$ so small that $\sum_{i=1}^l
\CW(\Sigma,B_\rho(p_i))\le \frac12 l \delta^2$ and hence we get
\begin{equation*}
  8\pi+\eps > \CW(\Sigma) +\tfrac12 l\delta^2 \ge 8\pi +\tfrac12 l\delta^2,
\end{equation*}
which is a contradiction for $\eps$ small enough.

Hence we can apply the approximate graphical decomposition lemma (in
the form of Lemma 3.4 in \cite{KMS:2011}) in order to conclude that
locally $\text{spt} \mu_k$ can be written as a multivalued graph away
from a small set of pimples.

The key result in order to get the regularity is a power-decay result
for the second fundamental form (see e.g. Lemma 3.6 in
\cite{KMS:2011}). Once we have this estimate, we can follow the rest
of the argument of \cite{KMS:2011} in order to conclude that
$\text{spt} \mu$ can locally be written as $C^{1,\alpha} \cap W^{2,2}$
graphs. Note that in our case we already ruled out the existence of
bad points and we also know that the limiting measure $\mu$ is coming
from the limiting immersion $F$. After having obtained this
preliminary regularity result, we can express $\CW$ in terms of the
graph functions and since $F$ is a minimizer of $\CW$ subject to the
area constraint, we conclude that graph functions solve the weak
Euler-Lagrange equation. Using the difference quotient technique as in
\cite{simon:1993existencewillmore} we finally get that $\text{spt}\mu$
and hence $F$ are smooth (in our case we get an additional lower order
term coming from the Lagrange parameter but this doesn't affect the
very general argument of Simon).

In order to get the power-decay result for $A_k$ we follow closely the
arguments of Lemma 3.6 in \cite{KMS:2011} and Lemma 5 in
\cite{Schygulla:2011}. More precisely we use the same replacement
procedures for the graph functions $u_k$ on balls of radii $\gamma \in
(\varrho/16,\varrho/32)$ as in the above mentioned lemmas in order to
get a comparison surface $\tilde{\Sigma}_k$ which satisfies
\begin{equation*}
  \big| |\tilde{\Sigma}_k|-a \big| \le c\varrho^2.
\end{equation*}
Using Lemma~\ref{thm:area_adjust_with_curvature} (note that $\Sigma_k
\subset B_{\rho/16}(p)$ and hence by construction we can assume that
$\tilde{\Sigma}_k \subset B_{\rho/4}(p)$) we conclude that there
exists a surface $\Sigma_k'$ with $|\Sigma_k'|=a$ and
\begin{equation*}
  \int_{\Sigma_k'} |A_k'|^2 \dmu_k'
  \le
  \int_{\tilde{\Sigma}_k} |\tilde{A}_k|^2 \tilde{\dmu}_k+C\varrho^2 a^{-1}.
\end{equation*}
These last two estimates allow us to use $\Sigma_k'$ as a comparison
sequence to the minimizing sequence $\Sigma_k$ and once we obtained
this fact we can follow the rest of the argument of Lemma 3.6 in
\cite{KMS:2011} word by word in order to get the desired power-decay.

Thus we have proved:
\begin{theorem}
  \label{thm:existence}
  Let $(M,g)$ be a compact, closed Riemannian manifold. Then there
  exists a constant $a_0>0$ such that for all $a\in(0,a_0)$ there is a
  smooth surface $\Sigma_a$ that minimizes the Willmore functional
  among all immersed surfaces with area $a$.
\end{theorem}


%
\section{The geometry of critical points}
\label{sec:geom-crit-points}
In this section we consider smooth solutions to the Euler-Lagrange
equation of the Willmore functional subject to an area constraint,
that is surfaces $\Sigma$ on which we have
\begin{equation}
  \label{eq:Wlambda}
  \Delta H + H |\Acirc|^2 + H \Ric(\nu,\nu) + H \lambda = 0.
\end{equation}
We show that if the area is small and the Willmore energy is close to
that of the Euclidean sphere, $\Sigma$ is very close to a geodesic
sphere. This can be used to conclude via the main result in
\cite{lamm-metzger:2010}, that $\Sigma$ is close to a critical point
of the scalar curvature.
\begin{proposition}
  \label{thm:curvest_l2_prop}
  Assume that $(M,g)$ has $C_B$-bounded geometry. Then there exist
  constants $C<\infty$ and $\eps>0$, depending only on $C_B$ such that
  the following holds.  If $\Sigma \subset (M,g)$ is a connected
  immersion and:
  \begin{enumerate}
  \item $\Sigma$ satisfies equation~\eqref{eq:Wlambda},
  \item $\lambda \leq \eps|\Sigma|^{-1}$,
  \item $\CW(\Sigma) \leq 8\pi + \eps$, and
  \item $|\Sigma|\leq \eps$.
  \end{enumerate}
  Then $\Sigma$ satisfies the following estimate:
  \begin{equation*}
    \int_\Sigma |\nabla^2 H|^2 + H^2 |\nabla H |^2 + H^4 |\Acirc|^2
    \dmu
    \leq
    C.
  \end{equation*}
\end{proposition}
For the proof we need the Bochner identity in the following form:
\begin{lemma}
  \label{thm:bochner}
  Let $(M,g)$ be a manifold and $\Sigma\subset M$ a smooth, compact,
  immersed 2-surface. Then for all $f\in C^\infty(\Sigma)$ we have
  that:
  \begin{equation*}
    \int_\Sigma |\nabla^2 f| ^2 \dmu
    =
    \int_\Sigma (\Delta f)^2 + |\nabla f|^2 \big( \half |\Acirc|^2 -
    \tfrac14 H^2 - \half \Scal + \Ric(\nu,\nu)\big) \dmu.
  \end{equation*}
  Here $\Scal$ and $\Ric$ denote the scalar and Ricci curvature of
  $(M,g)$.
\end{lemma}
\begin{proof}
  The Bochner identity states that
  \begin{equation*}
    \int_\Sigma |\nabla^2 f| ^2 \dmu
    =
    \int_\Sigma (\Delta f)^2 - \RicSig(\nabla f, \nabla f) \dmu    
  \end{equation*}
  where $\RicSig$ denotes the intrinsic Ricci curvature of $\Sigma$.
  Since $\RicSig(\nabla f, \nabla f) = \half \ScalSig |\nabla f|^2$ we
  can use the Gauss equation
  \begin{equation*}
    \half \ScalSig
    =
    \half  \Scal - \Ric(\nu,\nu) + \tfrac14 H^4 -
    \half |\Acirc|^2
  \end{equation*}
  to infer the claim.
\end{proof}
\begin{proof}[Proof of Proposition \ref{thm:curvest_l2_prop}]
  By the integrated Gauss equation we have that
  \begin{equation*}
    \CW(\Sigma)
    =
    8\pi
    + \int_\Sigma |\Acirc|^2 \dmu
    + 2\int_\Sigma G(\nu,\nu) \dmu
  \end{equation*}
  so that by the assumptions on the area and $\CW$ we obtain
  \begin{equation*}
    \|\Acirc\|_{L^2}^2 \leq C\eps.
  \end{equation*}
  Thus by choosing $\eps$ small we may later on assume that
  $\|\Acirc\|_{L^2}^2$ is as small as we desire.

  Furthermore, since $|\Sigma|$ is small, and $\CW$ is uniformly
  bounded, in view of lemma~\ref{thm:diam_leq_area}, we may also
  assume that $\Sigma\subset \CB_r(p)$ for some $p\in M$ and $r \leq C
  |\Sigma|^{1/2}$. Thus, we can use the Michael-Simon-Sobolev
  inequality as in \cite[Lemma 2.3]{lamm-metzger:2010} with a uniform
  constant on $\Sigma$, that is a constant that is at most double the
  one in Euclidean space, provided $\eps$ is small enough.
  
  Multiply equation~\eqref{eq:Wlambda} by $\Delta H$ and
  integrate. Integration by parts of the term including $\lambda$ and
  since $\Ric$ is bounded, we obtain:
  \begin{equation*}
    \int_\Sigma \half (\Delta H)^2 - \lambda |\nabla H|^2 \dmu
    \leq
    C \int_\Sigma H^2 |\Acirc|^4 + H^2 \dmu.
  \end{equation*}
  Thus
  \begin{equation}
    \label{eq:9}
    \int_\Sigma \half (\Delta H)^2\dmu
    \leq
    C + \eps|\Sigma|^{-1}\int_\Sigma |\nabla H |^2\dmu + C \int_\Sigma H^2 |\Acirc|^4\dmu.
  \end{equation}
  By the Bochner identity from Lemma \ref{thm:bochner} we infer the estimate
  \begin{equation}
    \label{eq:12}
    \int_\Sigma |\nabla^2 H|^2 + H^2 |\nabla H|^2 \dmu
    \leq
    C \int_\Sigma (\Delta H)^2 + (1 + |\Acirc|^2 ) |\nabla H|^2  \dmu.
  \end{equation}
  The Michael-Simon-Sobolev inequality implies that
  \begin{equation}
    \label{eq:10}
    \int_\Sigma |\nabla H|^2 \dmu
    \leq
    C|\Sigma| \int_\Sigma |\nabla^2 H |^2 + H^2|\nabla H|^2
      \dmu.    
  \end{equation}
  Inserting this and equation~\eqref{eq:9} into~\eqref{eq:12} yields
  \begin{equation*}
    \begin{split}
      &\int_\Sigma |\nabla^2 H|^2 + H^2 |\nabla H|^2 \dmu
      \\
      &\quad
      \leq
      C
      +
      C \int_\Sigma|\Acirc|^2 |\nabla H|^2 +H^2 |\Acirc|^4 \dmu
      \\
      &\qquad      
      +
      C (|\Sigma| + \eps) \int_\Sigma |\nabla^2 H |^2 + H^2|\nabla H|^2
        \dmu
    \end{split}
  \end{equation*}
  Thus, if $|\Sigma|$ and $\eps$ are small enough, we can absorb part
  of the right hand side of~\eqref{eq:12} to the left and infer
  together with~\eqref{eq:9} that
  \begin{equation}
    \label{eq:13}
    \int_\Sigma |\nabla^2 H|^2 + H^2 |\nabla H|^2 \dmu
    \leq
    C + C \int_\Sigma H^2 |\Acirc|^4 + |\Acirc|^2 |\nabla H|^2 \dmu.
  \end{equation}
  To proceed, recall the Simons identity, which implies that on an
  arbitrary immersed surface we have
  \begin{equation}
    \label{eq:simons}
    \begin{split}
      - \Acirc^{ij} \Delta \Acirc_{ij}
      + \half H^2 |\Acirc|^2
      &=
      -\la \Acirc, \nabla^2 H \ra
      +|\Acirc|^4
      +|\Acirc|^2 \Ric(\nu,\nu)
      \\
      & \phantom{=}
      - 2\Acirc^{ij}\Acirc_j^l \Ric_{il}
      - 2 \la \Acirc, \nabla\omega\ra.
    \end{split}      
  \end{equation}
  Here $\omega$ is the one-form $\omega=\Ric(\nu,\cdot)^T$ where the
  superscript $^T$ denotes projection to the tangential space of
  $\Sigma$. Multiply~\eqref{eq:simons} by $H^2$ and
  integrate. Integration by parts and the Codazzi equation $\div
  \Acirc = \half \nabla H + \omega$ yields
  \begin{equation}
    \label{eq:14}
    \begin{split}
      &\int_\Sigma H^2|\nabla \Acirc|^2 + 2 H \nabla _k H \Acirc^{ij}
      \nabla_k \Acirc_{ij} + \half H^4 |\Acirc|^2 \dmu
      \\
      &\quad
      \leq
      \int_\Sigma \la \div (H^2\Acirc), \nabla H + 2 \omega \ra \dmu
      + C + C \int_\Sigma H^2 |\Acirc|^4 \dmu.
    \end{split}
  \end{equation}
  Estimate
  \begin{equation}
    \label{eq:15}
    \left|\int_\Sigma 2 H \nabla _k H \Acirc^{ij} \nabla_k \Acirc_{ij}\dmu \right|
    \leq
    \frac{1}{4} \int_\Sigma H^2|\nabla \Acirc|^2 \dmu
    + 4 \int_\Sigma |\Acirc|^2 |\nabla H|^2 \dmu.
  \end{equation}
  Using the Codazzi equation to infer $\div \Acirc = \half \nabla H +
  \omega$ we calculate further
  \begin{equation}
    \label{eq:16}
    \begin{split}
      &\int_\Sigma \la \div (H^2 \Acirc), \nabla H + 2 \omega\ra \dmu
      \\
      &\quad
      =
      \int_\Sigma 2 H \Acirc (\nabla H, \nabla H + 2\omega) + H^2
      \la\div \Acirc , \nabla H + 2 \omega\ra \dmu
      \\
      &\quad
      =
      \int_\Sigma 2 H \Acirc(\nabla H, \nabla H + 2\omega)
      +\half H^2|\nabla H + 2\omega|^2\dmu
      \\
      &\quad
      \leq
      C + C \int_\Sigma H^2 | \nabla H |^2  + |\Acirc|^2 |\nabla H |^2\dmu.
    \end{split}
  \end{equation}
  Combining equations~\eqref{eq:14}, \eqref{eq:15} and \eqref{eq:16}
  yields
  \begin{equation}
    \label{eq:17}
    \int_\Sigma H^2 |\nabla A|^2 + H^4 |\Acirc|^2 \dmu
    \leq
    C + C \int_\Sigma H^2 |\nabla H |^2 + |\Acirc|^2 |\nabla H |^2 +
    H^2 |\Acirc |^4 \dmu.
  \end{equation}
  In view of~\eqref{eq:13} we finally infer
  \begin{equation}
    \label{eq:11}
    \begin{split}
      &\int_\Sigma |\nabla^2 H|^2  + |A|^2 |\nabla H|^2 +H^2|\nabla A|^2+ |A|^4 |\Acirc|^2
      \dmu
      \\
      &\quad
      \leq
      C + C \int_\Sigma H^2 |\Acirc|^4 + |\Acirc|^2|\nabla H|^2\dmu.
    \end{split}
  \end{equation}
  To proceed, we apply the Michael-Simon-Sobolev inequality and estimate
  \begin{equation}
    \label{eq:18}
    \begin{split}
      \int_\Sigma H^2|\Acirc|^4 \dmu
      &\leq
      C \left(\int_\Sigma |\nabla H| |\Acirc|^2 + |H||\Acirc||\nabla
        \Acirc| + H^2 |\Acirc|^2\dmu\right)^2
      \\
      &
      \leq
      \left(\int_\Sigma |\Acirc|^2\dmu\right)
      \left(\int_\Sigma |\Acirc|^2 |\nabla H|^2 + H^2|\nabla A|^2 +
        H^4|\Acirc|^2\dmu \right).
    \end{split}
  \end{equation}
  This shows that if $\eps$ and hence $\|\Acirc\|_{L^2}^2 $ is small,
  the first term on the right of equation~\eqref{eq:11} can be
  absorbed to the left.

  The second term on the right of~\eqref{eq:11} requires a little more
  work. By the Michael-Simon-Sobolev inequality we have
  \begin{equation}
    \label{eq:8}
    \begin{split}
      &\int_\Sigma | \Acirc |^2 |\nabla H|^2\dmu
      \\
      &\quad
      \leq
      C \int_\Sigma |\nabla \Acirc||\nabla H | + |\Acirc| |\nabla^2 H|
      + |H||\Acirc||\nabla H| \dmu
      \\
      &\quad
      \leq
      C \left( \int_\Sigma |\Acirc|^2\dmu \right)
      \left( \int_\Sigma |\nabla^2 H|^2 + H^2|\nabla H|^2 \dmu \right)
      + C \int_\Sigma |\nabla A|^2\dmu.
    \end{split}
  \end{equation}
  The first term on the right is of the same type as before. To
  estimate the second integrate the Simons identity~\eqref{eq:simons}
  to obtain
  \begin{equation}
    \label{eq:19}
    \int_\Sigma |\nabla \Acirc|^2 + \half H^2|\Acirc|^2 \dmu
    \leq
    C|\Sigma| + \int_\Sigma - \la \Acirc, \nabla^2 H\ra - 2\la \Acirc, \nabla
    \omega\ra + C |\Acirc|^4 \dmu.
  \end{equation}
  The Michael-Simon-Sobolev inequality implies that
  \begin{equation*}
    \int_\Sigma |\Acirc|^4\dmu
    \leq
    \left( \int_\Sigma |\Acirc|^2\dmu \right)
    \left( \int_\Sigma |\nabla \Acirc|^2 + H^2|\Acirc|^2 \dmu \right) .
  \end{equation*}
  so that the last term on the right of~\eqref{eq:19} can be absorbed
  to the left. Calculate further that
  \begin{equation*}
    -2 \int_\Sigma \la \Acirc, \nabla \omega\ra\dmu
    =
    2 \int_\Sigma \la \div \Acirc, \omega\ra \dmu
  \end{equation*}
  so that by Young's inequality
  \begin{equation*}
    \left| 2 \int_\Sigma \la \Acirc, \nabla \omega\ra\dmu\right|
    \leq
    \half \int_\Sigma |\nabla \Acirc|^2 \dmu
    + C \int_\Sigma |\omega|^2 \dmu.
  \end{equation*}
  Thus the second term on the right of equation~\eqref{eq:19} can also
  be absorbed to the left. Using Hölder's inequality we infer the estimate
  \begin{equation*}
    \int_\Sigma |\nabla \Acirc|^2 \dmu
    \leq
    C|\Sigma|
    +
    C \left( \int_\Sigma |\Acirc|^2 \dmu \right)
    \left( \int_\Sigma |\nabla^2 H|^2 \dmu \right).
  \end{equation*}
  Since the Codazzi-equation implies $\nabla H = 2\div \Acirc -
  2\omega$ we furthermore get
  \begin{equation*}
    \int_\Sigma |\nabla H|^2 \dmu
    \leq
    C|\Sigma| + C \int_\Sigma |\nabla \Acirc|^2 \dmu.
  \end{equation*}
  Thus the previous equation implies that
  \begin{equation}
    \label{eq:20}
    \int_\Sigma |\nabla A|^2\dmu
    \leq
    C|\Sigma|
    +
    C \left( \int_\Sigma |\Acirc|^2 \dmu \right)
    \left( \int_\Sigma |\nabla^2 H|^2 \dmu \right).
  \end{equation}
  Substituting this into equation~\ref{eq:8} yields that
  \begin{equation*}
    \begin{split}
      &\int_\Sigma |\Acirc|^2 |\nabla H|^2 \dmu
      \\
      &\quad
      \leq
      C|\Sigma| + 
      C \left( |\Sigma| + \int_\Sigma |\Acirc|^2 \dmu \right)
      \left( \int_\Sigma |\nabla^2 H|^2 + H^2|\nabla H|^2 \dmu \right).
    \end{split}
  \end{equation*}
  Thus we have shown that all but the constant term on the right of
  equation~\eqref{eq:13} can be absorbed to the left, provided $\eps$
  is small enough.
\end{proof}
We wish to complement these estimates by an estimate for
$\lambda$.
\begin{proposition}
  \label{thm:lambda_est}
  Let $(M,g)$ be a Riemannian manifold with $C_B$-bounded
  geometry. Let $\Sigma$ be a surface satisfying~\eqref{eq:Wlambda}
  for some $\lambda\in \IR$. Let $\rho_0$ be as in
  Remark~\ref{rem:uniform_normal_coordinates} and assume that $\Sigma
  \subset \CB_r(p)$ for some $p\in M$ and $0 < r < \rho_0$. Then we have
  the estimate
  \begin{equation*}
    |\lambda|
    \leq
    C|\Sigma|^{-1} \left( |\Sigma|^{1/2} + r \int_\Sigma |A|^2 \dmu \right).
  \end{equation*}
\end{proposition}
\begin{proof}
  The proof is based on the fact that the Willmore functional is scale
  invariant with respect to the Euclidean metric. Note that if
  $\Sigma$ is of Willmore type we have for all variations with normal
  velocity $f$ that
  \begin{equation*}
    \delta_f \CW(\Sigma) = \lambda \delta_f \CA(\Sigma)
  \end{equation*}
  where $\CA$ denotes the area functional. This implies that if
  $\delta_f\CA \neq 0$, we can write
  \begin{equation*}
    \lambda = \delta_f \CW / \delta_f \CA.
  \end{equation*}
  In Euclidean space we can choose $f=\la x, \nu \ra$ to be the normal
  velocity corresponding to scaling and infer that $\delta_f\CW =0$
  whereas $\delta_f \CA = 2 |\Sigma|$ so that in combination we get
  that $\lambda = 0$.

  In the situation as in the statement, this reasoning still works
  altough with some error terms.  Let $\Sigma \subset \CB_r(p)$ as in
  the statement of the proposition. Let $x$ denote the position vector
  field on $\CB_{\rho_0}(p)$ with respect to normal coordinates on
  $\CB_{\rho_0}(p)$. Since $\nabla_{\del_i} x^j = \delta_i^j +
  \Gamma_{ik}^jx^k$ it follows that
  \begin{equation}
    \label{eq:24}
    | \nabla x - \Id | \leq C r^2
    \qquad\text{and}\qquad
    |\nabla^2 x | \leq C.
  \end{equation}
  The first step is to estimate the variation of area with respect to
  a normal variation corresponding to scaling in normal
  coordinates. That is, we use $f=\la x,\nu\ra$ as a normal variation
  for $\Sigma$. This yields
  \begin{equation*}
    \delta_f \CA(\Sigma)
    =
    \int_\Sigma H \la x, \nu \ra \dmu
    =
    \int_\Sigma \div_\Sigma x \dmu.
  \end{equation*}
  Since $|\div_\Sigma x - 2| \leq C|x|^2$ by~\eqref{eq:24} and the
  fact that $\Sigma\subset \CB_r(p)$, we thus find that, provided $r$ is
  small enough,
  \begin{equation}
    \label{eq:21}
    \delta_{\la x,\nu\ra} \CA \geq |\Sigma|.
  \end{equation}
  As in the proof of Lemma~\ref{thm:curvature_scaling} we can
  estimate that
  \begin{equation*}
    |\delta_{\la x,\nu\ra} \CW|
    \leq
    C |\Sigma|^{1/2} \left( \int_\Sigma H^2 \dmu \right)^{1/2}
    +
    C r \int_\Sigma |A|^2 \dmu.
  \end{equation*}
  In combination with equation~\eqref{eq:21}, this yields the claimed
  estimate for~$\lambda$.  
\end{proof}
The main result of this section is a straight forward consequence of
the combination of propositions~\ref{thm:curvest_l2_prop}
and~\ref{thm:lambda_est}:
\begin{theorem}
  \label{thm:curvest_l2}
  Given a Riemannian manifold $(M,g)$ with $C_B$-bounded geometry
  there exist constants $C<\infty$ and $\eps>0$ depending only on
  $C_B$ such that the following holds.
  
  Assume that $\Sigma \subset (M,g)$ is a connected immersion that
  satisfies the following conditions:
  \begin{enumerate}
  \item $\Sigma$ satisfies equation~\eqref{eq:Wlambda},
  \item $\CW(\Sigma) \leq 8\pi + \eps$, and
  \item $|\Sigma|\leq \eps$.
  \end{enumerate}
  Then $\Sigma$ satisfies the following estimate:
  \begin{equation*}
    \int_\Sigma |\nabla^2 H|^2 + H^2 |\nabla H |^2 + H^4 |\Acirc|^2
    \dmu
    \leq
    C
  \end{equation*}
\end{theorem}
\begin{corollary}
  \label{thm:Hest_sup}
  Assume $\Sigma$ is as in theorem~\ref{thm:curvest_l2}. Then we have
  the following estimates:
  \begin{equation*}
    \|\Acirc\|_{L^2(\Sigma)} 
    \leq
    C|\Sigma|
  \qquad\text{and}\qquad
    \| H - 2/R \|_{L^\infty(\Sigma)}
    \leq
    C |\Sigma|^{1/2}
  \end{equation*}
  where $R$ is such that $|\Sigma|= 4\pi R^2$. In particular, if the
  area of $\Sigma$ is small enough we have that $H>0$.
\end{corollary}
\begin{proof}
  Apply the Michael-Simon-Sobolev inequality twice, first to estimate
  $H^2|\Acirc|^2$ and second to estimate $|\Acirc|^2$. This yields the
  first estimate.

  To see the second estimate, note that by the estimates of DeLellis
  and M\"uller we have that
  \begin{equation*}
    \| H - 2/R \|_{L^2(\Sigma)}
    \leq
    C \|\Acirc\|_{L^2(\Sigma)}
    \leq
    C |\Sigma|.
  \end{equation*}
  From \cite[Lemma 3.7]{lamm-metzger:2010} we infer that in addition
  \begin{equation*}
    \|H-2/R\|^4_{L^\infty(\Sigma)}
    \leq
    \|H- 2/R\|^2_{L^2(\Sigma)} \int_\Sigma |\nabla^2 H|^2 +
    H^4 |H-2/R|^2 \dmu
  \end{equation*}
  The first term in the integral is estimated by
  Theorem~\ref{thm:curvest_l2} and the second one can be absorbed to
  the left if $|\Sigma|$ is small enough. This yields the second
  estimate.
\end{proof}

\begin{corollary}
  Let $(M,g)$ be a compact, closed Riemannian manifold. Let $\Sigma_a$
  be the surfaces from Theorem \ref{thm:existence}. Then, if $a$ is
  small enough, $\Sigma_a$ has positive mean curvature.

  For any sequence $a_i \to 0$ there is a subsequence $a_{i'}$ such
  that $\Sigma_{a_{i'}}$ is asymptotic to a geodesic sphere centered
  at a point $p\in M$ where $\Scal$ attains its maximum.
\end{corollary}
\begin{proof}
  If $a$ is small enough, we find that Theorem~\ref{thm:curvest_l2}
  applies to $\Sigma_a$. Since these curvature estimates are all that
  is needed to carry out the analysis done in
  \cite{lamm-metzger:2010}, we obtain that $\Sigma_a$ is close to a
  geodesic sphere $S_r(p)$. From \cite[Theorem 5.1]{lamm-metzger:2010}
  we obtain the expansion
  \begin{equation*}
    \left|\CW(\Sigma_a) - 8\pi + \frac{|\Sigma_a|}{3}\Scal(p)\right|
    \leq
    Ca^{3/2}.
  \end{equation*}
  Comparing this expansion with the one for geodesic spheres,
  equation~\eqref{eq:27}, we find that the point $p$ is a global
  maximum of the scalar curvature of $(M,g)$.
\end{proof}


%
\appendix
\section{Complete surfaces of Willmore type}
\label{sec:classification}
In this section we use the methods developed in
\cite{fischer-schoen:1980} to classify complete surfaces of
Willmore type with positive mean curvature in Riemannian manifolds.

We start by recalling the Gauss equation
\begin{equation}
  \label{eq:c1}
  \ScalSig = \Scal - 2\Ric(\nu,\nu) + \half H^2 - |\Acirc|^2
\end{equation}
and the Euler-Lagrange equation satisfied by surfaces of Willmore type
\begin{equation}
  \label{eq:c2}
  \Delta H + H |\Acirc|^2 + H \Ric(\nu,\nu) + \lambda H = 0.
\end{equation}
Here $\lambda\in\IR$ is the Lagrange multiplier. Letting $f\in
C^1_c(\Sigma)$ and multiplying \eqref{eq:c2} with $f^2 H^{-1}$ we get
after integrating by parts
\begin{equation*}
  \int_\Sigma f^2\big(|\Acirc|^2 + \Ric(\nu,\nu) + \lambda+|\nabla \log
  H|^2\big)\dmu
  =
  2 \int_\Sigma f\la\nabla f,  \nabla \log H\ra\dmu.
\end{equation*}
Using Young's inequality we conclude
\begin{equation*}
  \int_\Sigma f^2\big(|\Acirc|^2 + \Ric(\nu,\nu) + \lambda\big)\dmu
  \leq
  \int_\Sigma |\nabla f|^2 \dmu.
\end{equation*}
Replacing the Ricci curvature on the left hand side by inserting
\eqref{eq:c1} we finally get the following lemma.
\begin{lemma}
  \label{stability}
  Let $\Sigma$ be a surface of Willmore type with positive mean
  curvature. Then we have for any $f\in C^1_c(\Sigma)$
  \begin{equation}
    \label{eq:c3}
    \int_\Sigma f^2\big(\tfrac12 |\Acirc|^2 + \tfrac14 H^2 +\tfrac12
    \Scal-\tfrac12 \ScalSig +\lambda\big)\dmu
    \leq
    \int_\Sigma |\nabla f|^2 \dmu.
  \end{equation}
  In particular, if $\lambda \geq -\frac12 \Scal$ we have
  \begin{equation}
    \label{eq:c4}
    \int_\Sigma f^2\big(\tfrac12 |\Acirc|^2 + \tfrac14 H^2 -\tfrac12 \ScalSig
    \big)\dmu
    \leq
    \int_\Sigma |\nabla f|^2 \dmu.
  \end{equation}
\end{lemma}
These inequalities are similar to the stability inequality for minimal
surfaces. Indeed they allow us to classify surfaces of Willmore type
with positive mean curvature. We directly get the following corollary:
\begin{corollary}
  \label{classificationcompact}
  Let $\Sigma \subset M$ be a compact surface of Willmore type with
  positive mean curvature and let $\lambda \geq -\frac12 \Scal$. Then
  $\Sigma$ is a topological sphere.
\end{corollary}
\begin{proof}
  In this situation we can insert $f\equiv 1$ into \eqref{eq:c4} and
  with the help of the Gauss-Bonnet theorem we get
\begin{equation}
  \label{eq:c5}
  0
  \leq
  \int_\Sigma \tfrac12 |\Acirc|^2 + \tfrac14 H^2\dmu
  \le
  4\pi (1-q(\Sigma)),
\end{equation}
where $q(\Sigma)$ is the genus of $\Sigma$. But if $\Sigma$ is a torus
then we conclude from the above inequality that $H\equiv 0$ which
contradicts the assumption of the corollary. This finishes the proof.
\end{proof}  
\begin{lemma}
  \label{noncompact}
  Let $\Sigma \subset M$ be a non-compact, complete surface of
  Willmore type with positive mean curvature and let $\lambda \geq
  -\frac12 \Scal$. Then $\Sigma$ is either conformally equivalent to
  the plane or to a cylinder. In the latter case $\Sigma$ has infinite
  absolute total curvature.
\end{lemma}
\begin{proof}
  We follow the closely the proof of Theorem~3 in
  \cite{fischer-schoen:1980}. Assume that the universal covering space
  of $\Sigma$ is $B_1(0)$. Defining $q=\frac12 \ScalSig-\frac12
  |\Acirc|^2 - \frac14 H^2$ and using~\eqref{eq:c4} we see that we can
  apply Lemma~1 of \cite{fischer-schoen:1980} and we get a positive
  solution $g$ on $\Sigma$ of the equation
  \begin{equation*}
    \Delta g-\tfrac12 \ScalSig g +(\tfrac12 |\Acirc|^2 + \tfrac14 H^2)g=0.
  \end{equation*}
  This solution can be lifted to $B_1(0)$ and by Corollary~3 in
  \cite{fischer-schoen:1980} this is a contradiction since $ \frac12
  \ScalSig=K$ and $\frac12 |\Acirc|^2 + \frac14 H^2\geq 0$.

  Hence the covering space of $\Sigma$ is $\mathbb{C}$ and this shows
  that $\Sigma$ is either a plane or a cylinder. If $\Sigma$ is a
  cylinder with finite absolute total curvature, then we can continue
  arguing as in Theorem~3 of \cite{fischer-schoen:1980} and we conclude
  \begin{equation*}
    \int_\Sigma  |\Acirc|^2 + \tfrac12 H^2\dmu\leq \int_\Sigma  \ScalSig\dmu \leq 0.
  \end{equation*}
  This contradicts the assumption $H>0$ and finishes the proof of the
  Lemma.
\end{proof}
For $M=\mathbb{R}^3$ and $\lambda=0$ we have the following theorem.
\begin{theorem}
  \label{completewillmore}
  Let $\Sigma\subset \mathbb{R}^3$ be a complete Willmore surface with
  positive mean curvature. Then $\Sigma$ is a round sphere.
\end{theorem}
\begin{proof}
  Combining Corollary~\ref{classificationcompact} and
  Lemma~\ref{noncompact} we conclude that $\Sigma$ is either a
  topological sphere or its universal cover is $\IC$. If $\Sigma$
  is a topological sphere we conclude from~\eqref{eq:c5} that
  \begin{equation*}
    \int_\Sigma |\Acirc|^2+\tfrac12 H^2\dmu\leq 8\pi.
  \end{equation*}
  Since on the other hand $\int_\Sigma\tfrac12 H^2\dmu \geq 8\pi$ for
  all closed surfaces $\Sigma \subset \mathbb{R}^3$ we find that
  $\Sigma$ is umbilic and hence a round sphere.

  Next we rule out the case that the universal cover of $\Sigma$ is
  $\IC$. We the Gauss equation~\eqref{eq:c1} yields that
  \begin{equation*}
    \tfrac12 \ScalSig=\tfrac14 H^2 -\tfrac12 |\Acirc|^2.
  \end{equation*}
  Inserting this into \eqref{eq:c4} we have for every $f\in
  C^1_c(\Sigma)$
  \begin{equation*}
    \int_\Sigma f^2|\Acirc|^2\dmu \leq \int_\Sigma |\nabla f|^2 \dmu.
  \end{equation*}
  Hence, defining $q=-|\Acirc|^2$, we can apply Theorem~1 of
  \cite{fischer-schoen:1980} and get a positive solution $g$ of
  \begin{equation*}
    \Delta g+ |\Acirc|^2g=0
  \end{equation*}
  on $\Sigma$ and by lifting also on $\mathbb{C}$. Hence we have a
  positive super-harmonic function $g$ on $\mathbb{C}$ which must be
  constant. This implies that $\Acirc\equiv 0$ and therefore $\Sigma$
  is a flat plane, which contradicts our assumptions.
\end{proof}

%
%
\bibliographystyle{abbrv}
\bibliography{../../extern/references}
\end{document}